\documentclass[a4,12pt]{article}
\begin{document}
\title{The Lindel{\"o}f Hypothesis for almost all 
Hurwitz's Zeta-Functions holds true                     
          }
\author{Masumi Nakajima \  \\
        \it Department of Economics \  \\
        \it International University of Kagoshima \   \\
        \it Kagoshima 891-0191, JAPAN  \\
        e-mail: nakajima@eco.iuk.ac.jp }
\maketitle
\begin{abstract}
By probability theory we prove here that the Lindel{\"o}f hypothesis holds for almost all Hurwitz's zeta-functions,
i.e. \\
$ \qquad \zeta({1\over2} + it,\omega)={\rm o}_{\omega,\epsilon}\{(\log t)^{{3\over2} + \epsilon}\} $ \\
for almost everywhere $ 0< \omega <1,$ and for any small $ \epsilon >0,$ where
 $ {\rm o}_{\omega,\epsilon} $ denotes the Landau small o-symbol which depends on
 $\omega$ and $\epsilon$ and $\zeta(s,\omega)$ denotes the Hurwitz zeta-function. 
The details will be given elsewhere.\\

Key words ; The Riemann zeta function, the Hurwitz zeta function, the Lindel{\"o}f hypothesis, 
            law of large numbers, law of the iterated logarithm.

Mathematics Subject Classification ; \\ 11M06, 11M26, 11M35, 60F15. 
\end{abstract}
Let $ \zeta(s,\omega)$ be the Hurwitz zeta function which is 
meromorphically  extended to the whole complex plane 
from the Dirichlet series 
\[ \sum_{n=0}^{\infty} (n+{\omega})^{-s} \quad 
(s=\sigma + it,\sigma=\Re s >1, \ 0<\omega \leq 1). \]
We should note that 
\[ \zeta(s,1)=\zeta(s),\]
\[ \zeta(s,{1\over2} )=({2^s}-1)\zeta(s), \]
where $\zeta(s)$ denotes the Riemann zeta function.\\
\quad In analytic number theory, there are three famous conjectures which are related each other as follows.\\
{\bf The Riemann Hypothesis }(1859, by B.Riemann):\\
$ \rho \notin {\bf R}, \ \zeta(\rho)=0 \Rightarrow \Re \rho ={1\over2} $  \\
{\bf The Lindel{\"o}f Hypothesis }(1908, by E.Lindel{\"o}f):\\
$ \zeta({1\over2}+ it)={\rm O}_{\epsilon}(t^{\epsilon}) \ for \ any \ small \ {\epsilon}>0,$ \\
where ${\rm O}_{\epsilon}$ denotes the Bachmann-Landau large O-symbol which depends on $\epsilon$. \\
{\bf The Density Hypothesis }: \\
\[ N(\sigma,T)={\rm O}_{\epsilon}(T^{2-2\sigma + \epsilon}) \]
\[ for \ any \ small \ \epsilon>0 \ and \ {1\over2}\leq \sigma \leq 1,\] \\
where $ N(\sigma,T) $ denotes the number of zeros of $\zeta(s)$ in the rectangle whose four vertices are $ \sigma ,1,1+iT $ and $ \sigma +iT $. \\
It is well known that \\
$the \ Riemann \ Hypothesis \Rightarrow the \ Lindel \ddot{o}f \ Hypothesis $ \\
$ \Rightarrow the \ Density \ Hypothesis.$ \\
(It is not known whether the Lindel{\"o}f Hypothesis implies 
the Riemann Hypothesis or not.)\\
And also as it is well known, the Riemann Hypothesis is the most important and the strongest conjecture that has serious influences on 
many branches of mathematics including number theory. 
But it is less known that in fact the Lindel{\"o}f Hypothesis has 
almost the same effects on number theory as the Riemann Hypothesis
~\cite{A1976}~\cite{I1985}~\cite{T1951}. About the Lindel{\"o}f 
Hypothesis there are many studies which improve the power $L$ of $t$
 in $\zeta({1\over2}+ it)={\rm O}(t^L)$. These studies in this 
direction have their long history and story. The recent results in 
this direction are due to G.Kolesnik, E.Bombieri, H.Iwaniec, 
M.N.Huxley, N.Watt and others, for example, 
$\zeta({1\over2}+ it)={\rm O}(t^{9/56})$ due to Bombieri and 
Iwaniec in 1986 and the best up to the present time is $={\rm O}(t^{32/205})$ due to 
Huxley in 2005. 
\\
\quad In 1952, Koksma and Lekkerkerker~\cite{K-L1952} proved that \\
\[ \int_0^1 |\zeta_{1}({1\over2}+it,\omega)|^{2}d\omega={\rm O}(\log t) \]
where $\zeta_{1}(s,\omega):=\zeta(s,\omega)-\omega^{-s}$ whose term $-\omega^{-s}$ makes keeping out the singularity at $ \omega=0$.
\\
\quad From this mean value results, by using {\v C}eby{\v s}ev's inequality in probability theory, 
we easily have 
\[ \mu \{ 0<\omega \leq 1 ; |\zeta({1\over2}+it,\omega)| \geq C \sqrt{\log t} \} \leq {{{\rm O}(1)}\over{C^2}} \]
\[ for \ any \ t>1 \ and \ any \ large \ C>0, \]
where $\mu \{ B \} $ denotes the Lebesgue measure of measurable set $ B $, which shows that the Lindel{\"o}f Hypothesis holds in the sence of weak law in probability theory. \\
\quad In this short note we give the following strong law version 
of the Lindel{\"o}f Hypothesis , that is, \\

\newtheorem{theo}{Theorem}
%
\begin{theo}

\begin{eqnarray*}
\zeta({1\over2}+it,\omega)={\rm o}_{\omega,\epsilon}\{(\log t)^{{3\over2} + \epsilon}\} 
\end{eqnarray*}
for almost everywhere  $ \omega \in \Omega:=(0,1) $ and for any small $ \epsilon >0 $.

\end{theo}

\quad In order to prove this theorem, we need some definitions and some results in probability theory. \\
\quad Let $ (\Omega, {  F},{\rm P}) $ be some probability space, 
$X,Y,Z,\cdots $ be complex valued random variables on this space, $ {\rm E}[X] $ be the 
expectation value of the random variable $ X $ and $ {\rm V}[X]={\rm E}[|X-{\rm E}[X]|^2] $ 
be the variance of $ X $.

\newtheorem{lem}{Lemma}
%
\begin{lem}
Let $ Z $ be a complex valued random variable. If $ {\rm E}[|Z|^2]|  
< +\infty $, then we have 
\[ |Z|<+\infty \ almost \ surely \ (abbreivated \ by \ a.s. ),\]
\[i.e. \ {\rm P}\{ |Z|<+\infty \}=1. \]
\end{lem}
{\bf proof.}\ From $|Z| \geq 0$, we have 
\[ 0 \leq |Z|=|Z(\omega)| < +\infty \ {\rm or} \ |Z|=|Z(\omega)| = +\infty. \]
We define the set $A \subset \Omega$ by 
$A:=\{ \omega ; |Z(\omega)|=+\infty \}$ 
and the indicator function of the set $A$;
\begin{eqnarray*}
1_{A}:=1_{A}(\omega):=\left\{
\begin{array}{ll}
1 & (\omega \in A)  \\
0 & (\omega \notin A).
\end{array}
\right.
\end{eqnarray*} 
If we assumed that ${\rm P}\{ \omega \in \Omega ; |Z(\omega)| < +\infty \} < 1 $, 
we would have ${\rm P}\{ A \}>0 $ and 
\[
{\rm E}[|Z|^2] \geq {\rm E}[|Z|^{2} 1_{A} ]=(+\infty){\rm P}\{ A \}=+\infty,
\]
which is the contradiction to the assumption $ {\rm E}[|Z|^2]|  
< +\infty. $ So we have the lemma.

%
\begin{lem}
Let $ Z_n $ be a complex valued random variables \ $(n=1,2,3.\cdots)$.
 If $ \sum_{n=1}^{\infty}{\rm E}[|Z_n|^2]  
< +\infty $, then we have \\
\[ Z_n \rightarrow 0 \ a.s.\ (as \ n \rightarrow +\infty ).\]
\end{lem}
{\bf proof.}\ By Lemma 1, we have 
\[
{\rm P}\{ \sum_{n=1}^{\infty}|Z_n|^2 < \infty \}=1,  
\]
which shows that $|Z_n|^2 \rightarrow 0 \ a.s.\ (as \ n \rightarrow +\infty )$, 
that is, $Z_n \rightarrow 0 \ a.s.\ (as \ n \rightarrow +\infty )$.
%
\begin{lem}
{\rm ({\rm Rademacher-Menchoff's lemma}~\cite{D1953}
~\cite{K1999}) } \\
Let $ a(p) \ (p=1,2,\cdots,2^{n+1}-1) $ be complex numbers and $ a(0):=0 $,  
then we have 
\begin{eqnarray*}
\lefteqn{ \max_{ 1\leq p <2^{n+1} } |a(p)|^2   }  \\
& \leq & (n+1) 
\sum_{k=0}^{n} \sum_{j=0}^{2^{n-k}-1}|a(2^{k}+j2^{k+1}) - a(j2^{k+1})|^2. 
\end{eqnarray*}
\end{lem}
{\bf proof.}\ For the natural number $p$ which satisfy $1 \leq p < 2^{n+1}$, 
we have its binomial expansion;
\[
p=\sum_{j=0}^{n}\epsilon_{j}2^{j} \ (\epsilon_{j}=0\ {\rm or} \ 1).
\]
With respect to the above $p$, we define $p_{k+1},\ p_{n+1},\ p_{0}$ 
respectively by
\begin{eqnarray*}
p_{k+1}:=\sum_{j=k+1}^{n}\epsilon_{j}2^{j} \ (k=0,1,2,\cdots,n-1),\\ 
p_{n+1}:=0,\ p_0:=p.
\end{eqnarray*}
From these definitions we have
\begin{eqnarray*}
p_0=p \geq p_1 \geq p_2 \geq \cdots \geq p_n \geq p_{n+1}=0, \ \ \ \ \ (1)
\\
 p_k - p_{k+1}=\epsilon_{k}2^k, \quad \quad \quad \\
 p_{k+1}=\sum_{j=k+1}^{n}\epsilon_{j}2^{j}
=\sum_{j=0}^{n-k-1}\epsilon_{k+1+j}2^{j+k+1} 
\\
=\sum_{j=0}^{n-k-1}(\epsilon_{k+1+j}2^{j})2^{k+1}=:
\delta_{k+1}2^{k+1},\ \ \ (2)
\\
0 \leq \delta_{k+1} \leq \sum_{i=0}^{n-k-1} 2^i = 2^{n-k}-1. 
\ \ \ (3)
\end{eqnarray*}
From 
\[
a(p)=a(p_0)=a(p_0)-a(p_{n+1})=\sum_{k=0}^{n}(a(p_{k})-a(p_{k+1})),
\]
we have
\begin{eqnarray*}
|a(p)|^2=|\sum_{k=0}^n 1 \cdot (a(p_k)-a(p_{k+1}))|^2 \\
\leq
\sum_{k=0}^n 1 \sum_{k=0}^n |a(p_k)-a(p_{k+1})|^2 \\
=(n+1)\sum_{k=0}^n |a(p_k)-a(p_{k+1})|^2 \\
=(n+1)\sum_{k=0}^n |a(\epsilon_k 2^k + p_{k+1})-a(p_{k+1})|^2 \ 
\\
 \ ({\rm by \ (1)}) \\
=(n+1)\sum_{k=0}^n |a(\epsilon_k 2^k + \delta_{k+1} 2^{k+1})-a(\delta_{k+1} 2^{k+1})|^2 \ 
\\
 \ ({\rm by \ (2)}) \\
\leq (n+1)\sum_{k=0}^n \sum_{j=0}^{2^{n-k}-1}
|a(2^k + j 2^{k+1})-a(j 2^{k+1})|^2, \  (4)
\\
\end{eqnarray*}
because we take the summation with respect to $k$ into account only 
when $\epsilon_k =1$, and we sum up $j$ in place of $\delta_k$ 
by (3). 
By the fact that the right hand side of (4) is independent of $p$, 
we have the lemma.

\newtheorem{de}{Definition}
\begin{de}
{\rm Let $ X,Y $ be complex valued random variables which satisfy \ 
$ {\rm E}[|X|^2],\ {\rm E}[|Y|^2]< \infty. $ \\
If $ {\rm E}[\bar{X}Y]={\rm E}[\bar{X}]{\rm E}[Y], $
\ we call $ X,Y $ (pairwise) uncorrelated.
}
\end{de}

\begin{de}
{\rm
Let $ X,Y $ be complex valued random variables which satisfy \ 
$ {\rm E}[|X|^2],\ {\rm E}[|Y|^2]< \infty. $ \\
If $ {\rm E}[\bar{X}Y]=0, $
\ we call $ X,Y $ {\rm(}pairwise{\rm)} orthogonal.
}
\end{de}

\begin{lem}
{\rm ({\rm Rademacher-Menchoff}~\cite{D1953}
~\cite{K1999}) } \\
Let $ X_1,X_2,\cdots $ be pairwise uncorrelated complex valued 
random variables which satisfy 
\[ {\rm E}[X_i]=0, \ \sigma_i^2:={\rm V}[X_i] \ (i=1,2,\cdots),
\ \sigma_i \geq 0 \]
and let $ S_n:=S_n(\omega):=X_1+X_2+\cdots+X_n. $ \\
Then we have 
\begin{eqnarray*}
\lefteqn{ {\rm E}[\max_{ 2^m < k \leq 2^{m+1} } |S_k-S_{2^m}|^2 ]   }  \\
& \leq & (m^2 + 1) 
\sum_{i=1}^{2^m} \sigma_{2^m + i}^2 \ for \ m=0,1,2,\cdots . 
\end{eqnarray*}
\end{lem}
{\bf proof.}\ In Lemma 3, we put 
\[
n=m-1, \  a(p)=X_{2^m+1}+X_{2^m+2}+\cdots+X_{2^m+p}.
\]
From this, we have 
\begin{eqnarray*}
 {\rm E}[\max_{ 2^m < k \leq 2^{m+1} } |S_k-S_{2^m}|^2 ]
 \\
\leq
{\rm E}[\max_{ 1 \leq k <2^{m} } |S_{2^m+k}-S_{2^m}|^2 ] 
+ {\rm E}[|S_{2^{m+1}}-S_{2^m}|^2 ] \\
={\rm E}[\max_{ 1 \leq p <2^{m} } |X_{2^m+1}+X_{2^m+2}+\cdots+X_{2^m+p}|^2 ] \\
+ {\rm E}[|X_{2^m+1}+X_{2^m+2}+\cdots+X_{2^{m+1}}|^2 ] \\
\leq 
{\rm E}[m\sum_{k=0}^{m-1}\sum_{j=0}^{2^{m-1-k}-1} |(X_{2^m+1}+
X_{2^m+2}+\cdots+X_{2^m+2^k+j2^{k+1}}) \\
-(X_{2^m+1}+X_{2^m+2}+\cdots+X_{2^m+j2^{k+1}})|^2 ] \\
+ {\rm E}[|X_{2^m+1}+X_{2^m+2}+\cdots+X_{2^{m+1}}|^2 ] \\
({\rm by \ Lemma \ 3})\\
=m \sum_{k=0}^{m-1}\sum_{j=0}^{2^{m-1-k}-1}
 {\rm E}[|X_{2^m+j2^{k+1}+1}+X_{2^m+j2^{k+1}+2}+ \\
\cdots
+X_{2^m+j2^{k+1}+2^k}|^2 ] \\
+ {\rm E}[|X_{2^m+1}+X_{2^m+2}+\cdots+X_{2^{m+1}}|^2 ]  
\\
=m \sum_{k=0}^{m-1}\sum_{j=0}^{2^{m-1-k}-1}
 {\rm V}[X_{2^m+j2^{k+1}+1}+X_{2^m+j2^{k+1}+2}+ \\
\cdots
+X_{2^m+j2^{k+1}+2^k} ] \\
+ {\rm V}[X_{2^m+1}+X_{2^m+2}+\cdots+X_{2^{m+1}} ] 
\\
=m \sum_{k=0}^{m-1}\sum_{j=0}^{2^{m-1-k}-1}
\sum_{i=1}^{2^k}\sigma_{2^m+j2^{k+1}+i}^2 
+\sum_{i=1}^{2^m}\sigma_{2^m+i}^2 \\
\leq m \sum_{k=0}^{m-1}
\sum_{i=1}^{2^m}\sigma_{2^m+i}^2 
+\sum_{i=1}^{2^m}\sigma_{2^m+i}^2
\\
\leq m \cdot m
\sum_{i=1}^{2^m}\sigma_{2^m+i}^2 
+\sum_{i=1}^{2^m}\sigma_{2^m+i}^2
\\
=(m^2+1)\sum_{i=1}^{2^m} \sigma_{2^m+i}^2,
\end{eqnarray*}
which completes the proof of the lemma.\\
\quad By using these lemmas, we have

\begin{theo}{\rm~\cite{N2004.1}} \\
Let $ X_1^{(n)},X_2^{(n)},\cdots , X_k^{(n)},\cdots $ be pairwise 
uncorrelated complex valued 
random variables which may depend on $ n $ and satisfy 
\[ {\rm E}[X_k^{(n)}]=0, \ \sigma_k^2:={\rm V}[X_k^{(n)}]=
{\rm O}(k^{-2\alpha}),\ |X_k^{(n)}|<+\infty \]
\[ (k,n=1,2,\cdots,\ \ \alpha \in {\bf R}, \forall \omega \in \Omega) \]
, where $ \sigma_k \geq 0 $ do not depend on $ n $ .
Also let 
\[ S_n^{(l)}:=S_n^{(l)}(\omega):= X_1^{(l)}+X_2^{(l)}+\cdots+X_n^{(l)}, \]
 and
\[ 
\varphi(n):=n^{\beta}(\log n)^{{3\over2}+\epsilon} 
\ with \ any \ small \ \epsilon >0 
\]
\begin{eqnarray*}
\beta :=\left\{
\begin{array}{ll}
0 & (\alpha \geq {1\over2})  \\
{1\over2}- \alpha & (\alpha < {1\over2}).
\end{array}
\right.
\end{eqnarray*}
Then we have
\[ S_n^{(n)}=S_n^{(n)}(\omega)= 
{\rm o}_{\omega,\epsilon}(\varphi(n)) \ a.s.\ \omega \in \Omega. \]
\end{theo}
{\bf proof.}\ We choose any natural number sequence $\{n_k\}_{k=1}^\infty$
 with $2^k<n_k\leq 2^{k+1}$ and 
$X_{1}^{(l)},X_{2}^{(l)},\cdots,X_{2^{m+1}}^{(l)},\cdots\ (l,m \in {\bf N})$
 are pairwise uncorrelated complex valued random variables for any 
$l \in {\bf N}$. We have 
\begin{eqnarray*}
{\rm E}[\sum_{k=1}^m |{{S_{2^k}^{(n_k)}}\over{\varphi(2^k)}}|^2
+\sum_{k=m+1}^\infty |{{S_{2^k}^{(2^{k+1})}}\over{\varphi(2^k)}}|^2]
=\sum_{k=1}^\infty 2^{-2k\beta}{(\log 2^k)}^{-3-2\epsilon}
(\sigma_1^2 + \sigma_2^2 + \cdots + \sigma_{2^k}^2). \ \ (5)\\
\end{eqnarray*}

In case of $\alpha > {1\over2}$, then $\beta=0$ and we have
\begin{eqnarray*}
(5)={\rm O}(\sum_{k=1}^\infty {(\log 2^k)}^{-3-2\epsilon}
\sum_{l=1}^{2^k} l^{-2\alpha}) \\
={\rm O}(\sum_{k=1}^\infty k^{-3-2\epsilon})<+\infty.
\end{eqnarray*}
In case of $\alpha = {1\over2}$, then $\beta=0$ and we have
\begin{eqnarray*}
(5)={\rm O}(\sum_{k=1}^\infty {(\log 2^k)}^{-3-2\epsilon}
\sum_{l=1}^{2^k} l^{-1}) \\
={\rm O}(\sum_{k=1}^\infty k^{-3-2\epsilon}\cdot \log 2^k) \\
={\rm O}(\sum_{k=1}^\infty k^{-2-2\epsilon})<+\infty.
\end{eqnarray*}
In case of $\alpha < {1\over2}$, then $\beta={1\over2}-\alpha$ and we have
\begin{eqnarray*}
(5)={\rm O}(\sum_{k=1}^\infty 2^{-k(1-2\alpha)} 
{(\log 2^k)}^{-3-2\epsilon}
\sum_{l=1}^{2^k} l^{-2\alpha}) \\
={\rm O}(\sum_{k=1}^\infty 2^{-k(1-2\alpha)} k^{-3-2\epsilon}
\cdot  2^{k(1-2\alpha)}) \\
={\rm O}(\sum_{k=1}^\infty k^{-3-2\epsilon})<+\infty.
\end{eqnarray*}
Then in any case, we have
\[
{\rm E}[\sum_{k=1}^m |{{S_{2^k}^{(n_k)}}\over{\varphi(2^k)}}|^2
+\sum_{k=m+1}^\infty |{{S_{2^k}^{(2^{k+1})}}\over{\varphi(2^k)}}|^2]
<+\infty ,
\]
which means, by Lemma 2, with some $A((n_1,\cdots,n_m))\subset\Omega$,
\[
\sum_{k=1}^m |{{S_{2^k}^{(n_k)}(\omega)}\over{\varphi(2^k)}}|^2
+\sum_{k=m+1}^\infty |{{S_{2^k}^{(2^{k+1})}(\omega)}\over{\varphi(2^k)}}|^2
<+\infty
\]
for $\forall \omega \in A((n_1,\cdots,n_m))$ with 
${\rm P}\{A((n_1,\cdots,n_m))\}=1$.
We put
\[
A(m):=\bigcap_{(n_1,\cdots,n_m)}A((n_1,\cdots,n_m))
\]
where $(n_1,\cdots,n_m)$ under $\cap$ runs through all 
$(n_1,\cdots,n_m)\in {\bf N}^{m}$ with 
$2^k<n_k\leq 2^{k+1} (k=1,2,\cdots,m).$ \\
Since 
\[\bigcap_{(n_1,\cdots,n_m)}\]
is finitely many intersections of the sets, we have
\[{\rm P}\{A(m)\}=1.\]
Therefore we have
\[
\sum_{k=1}^m |{{S_{2^k}^{(n_k)}(\omega)}\over{\varphi(2^k)}}|^2
+\sum_{k=m+1}^\infty |{{S_{2^k}^{(2^{k+1})}(\omega)}\over{\varphi(2^k)}}|^2
<+\infty
\]
for $\forall\omega \in A(m)$ with \ ${\rm P}\{A(m)\}=1$.\\
We show that
\[A(m)=A(m+1) \ (m=1,2,\cdots).\]
In fact, if $\omega \in A(m)$ which means
\[
\sum_{k=1}^m |{{S_{2^k}^{(n_k)}(\omega)}\over{\varphi(2^k)}}|^2
+\sum_{k=m+1}^\infty |{{S_{2^k}^{(2^{k+1})}(\omega)}\over{\varphi(2^k)}}|^2
<+\infty
\],
then we immediately have
\[
\sum_{k=1}^m |{{S_{2^k}^{(n_k)}(\omega)}\over{\varphi(2^k)}}|^2
+|{{S_{2^{m+1}}^{(n_{m+1})}(\omega)}\over{\varphi(2^{m+1})}}|^2
+\sum_{k=m+2}^\infty |{{S_{2^k}^{(2^{k+1})}(\omega)}\over{\varphi(2^k)}}|^2
<+\infty
\]
for $2^{m+1}<\forall n_{m+1}\leq 2^{m+2}$, 
because $|X_k^{(l)}(\omega)|<+\infty$ for $\forall \omega \in \Omega.$
This means $\omega \in A(m+1)$. Inversely $\omega \in A(m+1)$ implies 
$\omega \in A(m)$ by the same argument.\\
So, There exists
\[\lim_{m \to +\infty}A(m)=:A=A(1) \ \ {\rm and} \ \ {\rm P}\{A\}=1.\]
This means 
\[
\sum_{k=1}^\infty |{{S_{2^k}^{(n_k)}(\omega)}\over{\varphi(2^k)}}|^2
<+\infty \ for \ \forall\{n_k\}_{k=1}^\infty,\ \forall \omega \in A \ 
with \ {\rm P}\{A\}=1
\]
and
\[
\lim_{k \to +\infty}{{S_{2^k}^{(n_k)}(\omega)}\over{\varphi(2^k)}}
=0 \ a.s. \ \omega \ for \ \forall\{n_k\}_{k=1}^\infty \quad \quad 
\quad {\rm (6)}
\]
Next we put 
\begin{eqnarray*}
Y_k^{(n_k)}:=\max_{1\leq l \leq 2^k}
      |X_{2^k+1}^{(n_k)}+X_{2^k+2}^{(n_k)}+\cdots+X_{2^k+l}^{(n_k)}| \\
.
\end{eqnarray*}
By Lemma 4, we have for any $l \in {\bf N}$
\begin{eqnarray*}
{\rm E}[|Y_k^{(l)}|^2]\leq (k^2+1)\sum_{i=1}^{2^k}\sigma_{2^k+i}^2 \\
=\left\{
\begin{array}{ll}
{\rm O}(k^2+1) & (\alpha \geq {1\over2})  \\
{\rm O}(2^{(k+1)(1-2\alpha)}(k^2+1)) & (\alpha < {1\over2}).
\end{array}
\right.
\end{eqnarray*}
and
\[
{\rm E}[\sum_{k=1}^m |{{Y_{k}^{(n_k)}}\over{\varphi(2^k)}}|^2
+\sum_{k=m+1}^\infty |{{Y_{k}^{(2^{k+1})}}\over{\varphi(2^k)}}|^2]
={\rm O}(\sum_{k=1}^\infty 2^{-2k\beta} k^{-3-2\epsilon}
{\rm E}[ |Y_{k}^{(l)}|^2])\ \ with \ \forall l\in {\bf N}. \ \ (7)
\]
In case of $\alpha \geq {1\over2}$, then $\beta=0$ and we have
\begin{eqnarray*}
(7)={\rm O}(\sum_{k=1}^\infty k^{-3-2\epsilon}(k^2+1) )
={\rm O}(\sum_{k=1}^\infty k^{-1-2\epsilon})<+\infty.
\end{eqnarray*}

In case of $\alpha < {1\over2}$, then $\beta={1\over2}-\alpha$ and we have
\begin{eqnarray*}
(7)={\rm O}(\sum_{k=1}^\infty 2^{-k(1-2\alpha)} k^{-3-2\epsilon}
\cdot (k^2+1) 2^{(k+1)(1-2\alpha)}) \\
={\rm O}(2^{-2\alpha}\sum_{k=1}^\infty k^{-1-2\epsilon})<+\infty.
\end{eqnarray*}
In any case, we have
\[
{\rm E}[\sum_{k=1}^m |{{Y_{k}^{(n_k)}}\over{\varphi(2^k)}}|^2
+\sum_{k=m+1}^\infty |{{Y_{k}^{(2^{k+1})}}\over{\varphi(2^k)}}|^2]
<+\infty .
\]
The same argument as that of ${{S_{2^k}^{(n_k)}}\over{\varphi(2^k)}}$ 
leads
\[
\lim_{k\to\infty}{{Y_{k}^{(n_k)}(\omega)}\over{\varphi(2^k)}}
=0 \ {\rm a.s.} \ \ {\rm for}\ \forall\{n_k\}_{k=1}^\infty . \ \ \ (8)
\]
Then, for any $n$ with $2^m<n\leq 2^{m+1}$, we have, by (6) and (8),
\begin{eqnarray*}
{{|S_{n}^{(n)}|}\over{\varphi(n)}} \leq 
{{|S_{n}^{(n)}|}\over{\varphi(2^m)}} \leq 
{ { |S_{2^m}^{(n)}|+Y_{m}^{(n)} }\over{\varphi(2^m)} }\\
={ { |S_{2^m}^{(n)}| }\over{\varphi(2^m)} }+
 { { Y_{m}^{(n)} }\over{\varphi(2^m)} }
\to 0 \ \ {\rm a.s.}\ \ ({\rm as}\ n\to\infty),
\end{eqnarray*}
which means
\[ S_n^{(n)}=S_n^{(n)}(\omega)= 
{\rm o}_{\omega}(\varphi(n)) \ a.s.\ \omega \in \Omega, 
\ \ {\rm with \ any \ small \ \epsilon>0.}
\]
This completes the proof.

\newtheorem{rem}{Remark}
\begin{rem}
{\rm
 This theorem is a generalization of the strong 
limit theorem the position of which may be placed between laws of large numbers and 
laws of the iterated logarithm in probability theory. 
(\ Therefore we would like to call these types of theorems quasi laws 
of the iterated logarithm.)
}
\end{rem}

\begin{rem}
{\rm
 This is also a new proof of the strong law of large numbers 
without using the Borel-Cantelli theorem. We can prove other 
limit theorems in probability theory by this method.
}
\end{rem}

\  We yet need some lemmas for proving Theorem 1. 

\begin{lem}
{\rm ( {\rm Functional equation for the Hurwitz zeta function}
~\cite{A1976}~\cite{I1985}~\cite{T1951} ) } \\

\[ \zeta(s,\omega)={{\Gamma(s)}\over{(2\pi)^s}}
\{ e^{-{\pi\over2}is}{\rm F}(\omega,s) + 
e^{+{\pi\over2}is}{\rm F}(-\omega,s) \} \]
for \ $ 0<\omega<1,\sigma>0 $ or \ $ 0<\omega \leq 1,\sigma>1 $,\\
where $\Gamma(s)$ is the gamma function of Euler and \\
${\rm F}(\omega,s):=\sum_{k=1}^{\infty}k^{-s}e^{2\pi ik \omega}.$
\end{lem}

\begin{lem}
{\rm ~\cite{I1985}~\cite{M1970}~\cite{T1980}~\cite{T1951}}
\[ |\Gamma(s)|=\sqrt{2\pi}|t|^{\sigma -{1\over2}}e^{-{\pi\over2}|t|}
\{ 1 + {\rm O}_{\sigma_1,\sigma_2,\delta}({1\over{|t|}}) \} \]

for \ $ \sigma_1 \leq \sigma \leq \sigma_2,\ |t| \geq \delta >0. $

\end{lem}

\begin{lem}

\[ 
{\rm F}(\omega,s)=\sum_{k \leq t^2}k^{-s}e^{2\pi ik \omega} + 
{\rm O}({1\over{|1-e^{2\pi i \omega}|}}t^{1-2\sigma})
 \]
for \ $ \sigma \geq {1\over2} . $

\end{lem}
{\bf proof.}\ 
\[ 
{\rm F}(\omega,s)=\sum_{k \leq t^2}k^{-s}e^{2\pi ik \omega} + 
\sum_{k > t^2}k^{-s}e^{2\pi ik \omega}.
 \]
By applying the partial summation to the second term of the above 
${\rm F}(\omega,s)$,
\begin{eqnarray*}
\sum_{k > t^2}k^{-s}e^{2\pi ik \omega}\\
=-A(t^2){(t^2)}^{-s}+s\int_{t^2}^\infty
A(u)u^{-s-1}du \\
={\rm O}({1\over{|1-e^{2\pi i \omega}|}}t^{-2\sigma})+
{\rm O}(t{1\over{|1-e^{2\pi i \omega}|}}t^{-2\sigma}) \\
={\rm O}({{t^{1-2\sigma}}\over{|1-e^{2\pi i \omega}|}}),
\end{eqnarray*}
where $A(t^2):=\sum_{k < t^2}e^{2\pi ik \omega}$, which completes 
the proof of the lemma.

From Lemma 5,6 we easily have
\begin{lem}

\[ 
|\zeta({1\over2} + it,\omega)|={\rm O}({\rm F}
(\omega,{1\over2} + it))
\ \ for \ 0<\omega<1.
\]

\end{lem}

\noindent
{\bf Proof of Theorem 1.}
From Lemma 7, we have 
\[
{\rm F}(\omega,{1\over2}+it)={\rm O}_{\delta}
(\sum_{k \leq t^2}k^{-{1\over2}-it}e^{2\pi ik \omega}) 
\]
\[ for \ 0<\delta<\omega<1-\delta \ with \ any \ small \ \delta>0.
\] 
In Theorem 2,  put 
$ \Omega=(0,1),\ {\rm P}=\mu $ (Lebesgue \ measure),  
$ n=[t^2] $ \ ($ [x] $ denotes the integral part of 
real number $ x.$)  and 
\[
X_k^{(t)}=k^{-{1\over2}-it}e^{2\pi ik \omega} \ 
\ (k=1,2,\cdots), 
\]
which satisfy all the conditions in Theorem 2. Then we have
\[
\sum_{k \leq t^2}k^{-{1\over2}-it}e^{2\pi ik \omega}
={\rm o}_{\omega,\epsilon}(\varphi([t^2]))
={\rm o}_{\omega,\epsilon}((\log t)^{{3\over2}+\epsilon})
\]
With Lemma 7,8, this completes the proof of the theorem. 

\begin{rem}
{\rm
 The exact expression of the Lindel{\"o}f Hypothesis is 

\begin{eqnarray*}
\mu_{\omega}(\sigma)=\left\{
\begin{array}{ll}
0 & (\sigma \geq {1\over2})  \\
{1\over2}- \sigma & (\sigma < {1\over2}).
\end{array}
\right.
\end{eqnarray*}
\begin{eqnarray*}
where \ \ \mu_{\omega}(\sigma):=
\lim_{\stackrel{\scriptstyle }{t\ \to}}\sup_{+\infty}
{{\log |\zeta(\sigma + it,\omega)|}\over{\log t}}, 
\end{eqnarray*}
which is the same form as

\begin{eqnarray*}
\beta =\left\{
\begin{array}{ll}
0 & (\alpha \geq {1\over2})  \\
{1\over2}- \alpha & (\alpha < {1\over2}),
\end{array}
\right.
\end{eqnarray*}
in Theorem 2. \\
}
\end{rem}

\begin{rem}
{\rm
 In 1936, Davenport and Heilbronn~\cite{D-H1936} has already proved 
that the Riemann Hypothesis fails for $\zeta(s,\omega)$ with 
transcendental number $\omega$ and rational number $\omega \neq 
{1\over2},1$ in contrast with our Theorem 1, which shows that 
the Lindel{\"o}f Hypothesis by itself, for example, without 
the Euler product, does not imply the Riemannn 
Hypothesis.
}
\end{rem}
\begin{rem}
{\rm
 It seems that the behaviour of 
$\zeta(s,\omega)$ as $\omega$ varies in the interval $(0,1)$ 
is very complicated because of the following facts;\\
(1)Barasubramanian-Ramachandra~\cite{B-R1977}(\ the case $ \omega=1
 $ ) and 
Ramachandra-Sankaranarayanan~\cite{R-S1989} proved the following 
$\Omega$-theorem;
\begin{eqnarray*}
\zeta({1\over2} + it,\omega)=\Omega(\exp (C_{\omega}\sqrt{{\log t}
\over{\log\log t}})) \\ 
with \ some \ C_{\omega}>0 \ and \ \omega \in {\bf Q},
\end{eqnarray*}
which shows 
\begin{eqnarray*}
\{0<\omega<1;{\rm  Theorem \ 1. \  holds } \} 
\cap {\bf Q}=\emptyset.
\end{eqnarray*}
\noindent
(2)It is well known that divisor problems and circle problems 
are closely related each other and so are shifted divisor problems 
and shifted circle problems.
The Hurwitz zeta function naturally appears in shifted divisor 
problems~\cite{N1993}.
And Bleher-Cheng-Dyson-Lebowitz
~\cite{B-C-D-L1993} pointed 
out that the value distributions of the error terms of the number of 
lattice points inside shifted circles behave very differently when 
the shift varies by their numerical studies. 
Therefore it seems that the behaviour of $\zeta(s,\omega)$ including 
its value distribution is very complicated as $\omega$ varies.
(For the value distribution of $\zeta(s,\omega)$ with transcendental 
number $\omega$, see ~\cite{N1997}. )
\\
\noindent
(3)Our numerical studies by "Mathematica"  show also the 
complexity of the behaviour of $\zeta(s,\omega)$ as follows, 
for example,}
\end{rem}
\noindent
The graph of $\zeta(s,x)$ which plots the points
$(x,y) \in {\bf R^2}$ such that 
\[
y={{|\zeta({1\over2}+it,x)-x^{-({1\over2}+it)}|}\over{(\log t)^2}}
\ (0 \leq x \leq 1,\ t=10^8.).
\]
seems to be a kind of white noise.\\
\noindent
{\bf Acknowledgment} \ \ The author thanks to \\
Prof. Jyoichi Kaneko of the University of the Ryukyus for 
his careful reading of the previous manuscript.\\

\end{document}